\documentclass[10pt]{article}
\usepackage{epic,latexsym,amssymb}
\usepackage{color}
\usepackage{tikz}
\usepackage{amsfonts,epsf,amsmath}

\newtheorem{thm}{Theorem}
\newtheorem{lem}{Lemma}

\newcommand{\qed}{$\Box$}

\newcommand{\proof}{\noindent\textbf{Proof. }}

\let\oldenumerate\enumerate
\renewcommand{\enumerate}{
  \oldenumerate
  \setlength{\itemsep}{0pt}
  \setlength{\parskip}{0pt}
  \setlength{\parsep}{0pt}
}

\begin{document}

\title{Hamiltonicity of Domination Vertex-Critical Claw-Free Graphs}

\author{Pawaton Kaemawichanurat\thanks{This work was funded by Development and Promotion of Science Technology Talents (DPST) Research Grant 031/2559.}
\\ \\
Theoretical and Computational Science Center \\
and Department of Mathematics\\
King Mongkut's University of Technology Thonburi \\
Bangkok, Thailand \\
\small \tt Email: pawaton.kae@kmutt.ac.th}
\date{}

\maketitle

\begin{abstract}
A graph $G$ is said to be $k$-$\gamma$-vertex critical if the domination numbers $\gamma(G)$ of $G$ is $k$ and $\gamma(G - v) < k$ for any vertex $v$ of $G$. Similarly, A graph $G$ is said to be $k$-$\gamma_{c}$-vertex critical if the connected domination numbers $\gamma_{c}(G)$ of $G$ is $k$ and $\gamma_{c}(G - v) < k$ for any vertex $v$ of $G$. The problem of interest is to determine whether or not $2$-connected $k$-$\gamma$-vertex critical graphs are Hamiltonian. In this paper, for all $k \geq 3$, we provide a $2$-connected $k$-$\gamma$-vertex critical graph which is non-Hamiltonian. We prove that every $2$-connected $3$-$\gamma$-vertex critical claw-free graph is Hamiltonian and the condition claw-free is necessary. For $k$-$\gamma_{c}$-vertex critical graphs, we present a new method to prove that every $2$-connected $3$-$\gamma_{c}$-vertex critical claw-free graph is Hamiltonian. Moreover, for $4 \leq k \leq 5$, we prove that every $3$-connected $k$-$\gamma_{c}$-vertex critical claw-free graph is Hamiltonian. We show that the condition claw-free is necessary by giving $k$-$\gamma_{c}$-vertex critical non-Hamiltonian graphs containing a claw as an induced subgraph for $3 \leq k \leq 5$.
\end{abstract}

{\small \textbf{Keywords:}domination, critical, Hamiltonian.
\\
\indent {\small \textbf{AMS subject classification:} 05C69, 05C45

\section{\bf Introduction}
In this paper, let $G$ denote a finite graph with vertex set $V(G)$ and edge set $E(G)$. For $S \subseteq V(G)$, $G[S]$ denotes the subgraph of $G$ induced by $S$. Throughout this paper all graphs are simple (no loops or multiple edges) and connected. The \emph{neighborhood} $N_{G}(x)$ of a vertex $x$ in $G$ is the set of vertices of $G$ which are adjacent to $x$. For vertex subsets $X, Y \subseteq V(G)$, $N_{Y}(X)$ denotes the set of vertices $y \in Y$ which $yx \in E(G)$ for some $x \in X$. For a subgraph $H$ of $G$, we use $N_{Y}(H)$ instead of $N_{Y}(V(H))$ and we use $N_{H}(X)$ instead of $N_{V(H)}(X)$. If $X = \{x\}$, we use $N_{Y}(x)$ instead of $N_{Y}(\{x\})$. A vertex $v$ of $G$ is an \emph{isolated vertex} if $N_{G}(v) = \emptyset$. Moreover, if $N_{G}(v) = \{u\}$, then $v$ is called an \emph{end vertex} and $u$ is called a \emph{support vertex}. A \emph{tree} is a connected graph having no cycle. A \emph{star} $K_{1, n}$ is a tree of order $n + 1$ having $n$ end vertices. Let $u_{1}u_{2}u_{3}...u_{s_{1} + 1}, v_{1}v_{2}v_{3}...v_{s_{2} + 1}$ and $w_{1}w_{2}w_{3},...w_{s_{3} + 1}$ be three disjoint paths of length $s_{1}, s_{2}$ and $s_{3}$ respectively. A \emph{net} $N_{s_{1}, s_{2}, s_{3}}$ is constructed by adding edges so that the vertices in $\{u_{s_{1} + 1}, v_{s_{2} + 1},$ $w_{s_{3} + 1}\}$ form a cycle of length three. For a finite family $\mathcal{G}$ of graphs, a graph $G$ is $\mathcal{G}$\emph{-free} if $G$ does not contain $H$ as an induced subgraph for any graph $H$ in $\mathcal{G}$.
\vskip 5 pt

\indent An \emph{independent set} of $G$ is a set of pairwise non-adjacent vertices and the \emph{independence number} of $G$ is the maximum cardinality of an independent set. Let $\omega(G)$ be the number of components of $G$. For a connected graph $G$, a \emph{cut set} $S$ is a subset of $V(G)$ such that $\omega(G - S) > 1$, moreover, if $S = \{v\}$, then $v$ is called a \emph{cut vertex}. The \emph{connectivity} $\kappa(G)$ of $G$ is the minimum $s$ for which $G$ has a cut set of $s$ vertices. If $G$ has no cut set, then we define $\kappa(G)$ to be $|V(G)| - 1$. A graph $G$ is \emph{$l$-connected} if $l \leq \kappa(G)$. A \emph{Hamiltonian cycle(path)} is a cycle(path) containing every vertex of a graph. A graph $G$ is \emph{Hamiltonian} if it contains a Hamiltonian cycle.
\vskip 5 pt

\indent For subsets $D, X \subseteq V(G)$, $D$ \emph{dominates} $X$ if every vertex in $X$ is either in $D$ or adjacent to a vertex in $D$. If $D$ dominates $X$, then we write $D \succ X$, further, we write $a \succ X$ when $D = \{a\}$. If $X = V(G)$, then $D$ is a \emph{dominating set} of $G$ and we write $D \succ G$ instead of $D \succ V(G)$. A smallest dominating set is called a $\gamma$-\emph{set}. The \emph{domination number} $\gamma(G)$ of $G$ is the cardinality of a $\gamma$-set of $G$. A \emph{connected dominating set} of a graph $G$ is a dominating set $D$ of $G$ such that $G[D]$ is connected. If $D$ is a connected dominating set of $G$, we then write $D \succ_{c} G$. The minimum cardinality of a connected dominating set of $G$ is called the \emph{connected domination number} of $G$ and is denoted by $\gamma_{c}(G)$. A \emph{total dominating set} of a graph $G$ is a subset $D$ of $V(G)$ such that every vertex in $G$ is adjacent to a vertex in $D$. The minimum cardinality of a total dominating set of $G$ is called the \emph{total domination number} of $G$ and is denoted by $\gamma_{t}(G)$.
\vskip 5 pt

\indent A graph $G$ is said to be $k$-$\gamma$-\emph{edge critical} if $\gamma(G) = k$ and $\gamma(G + uv) < k$ for any pair of non adjacent vertices $u, v$ of $G$. We define a $k$-$\gamma_{c}$-\emph{edge critical graph} and a $k$-$\gamma_{t}$-\emph{edge critical} similarly.
\vskip 5 pt

\indent A graph $G$ is said to be $k$-$\gamma$-\emph{vertex critical} if $\gamma(G) = k$ and $\gamma(G - v) < k$ for every vertex $v$ of $G$. We define a $k$-$\gamma_{c}$-\emph{vertex critical graph} similarly. Obviously, a disconnected graph does not have a connected dominating set. If a graph has connectivity one, then it is disconnected after removing a cut vertex. Thus, when we study on $k$-$\gamma_{c}$-vertex critical graphs, we always focus on $2$-connected graphs. A graph $G$ is said to be $k$-$\gamma_{t}$-\emph{vertex critical} if $\gamma_{t}(G) = k$ and $\gamma_{t}(G - v) < k$ for any vertex $v$ of $G$ which is not a support vertex.
\vskip 5 pt

\indent The classical work on the Hamiltonicity of $k$-$\gamma$-edge critical graphs has been discussed since 1983 by Sumner and Blitch~\cite{SB}. Then, Wojcicka~\cite{W} proved that every connected $3$-$\gamma$-edge critical graph of at least $7$ vertices has a Hamiltonian path. Further, Flandrin et al.~\cite{FTWZ}, Favaron et al.~\cite{FTZ} together with Tian et al.~\cite{FWZ} proved that all connected $3$-$\gamma$-edge critical graph with $\delta \geq 2$ are Hamiltonian. In 1998, Sumner and Wojcicka (\cite{HHS} Chapter 16) conjectured that: ``For $k \geq 4$, all $(k - 1)$-connected $k$-$\gamma$-edge critical graphs are Hamiltonian." In 2005, Yuansheng et al.~\cite{YCXYX} disproved this conjecture by giving a $3$-connected $4$-$\gamma$-edge critical non-Hamiltonian graph. However, Kaemawichanurat and Caccetta~\cite{KC} proved that this conjecture is true if $k = 4$ and the graphs are claw-free. That is : every $3$-connected $4$-$\gamma$-edge critical claw-free graph is Hamiltonian.
\vskip 5 pt

\indent In the context of $k$-$\gamma_{c}$-edge critical graphs, Kaemawichanurat et al.~\cite{KCA2} showed that, for $1 \leq k \leq 3$, every $2$-connected $k$-$\gamma_{c}$-edge critical graph is Hamiltonian. They, further, gave a construction of $k$-$\gamma_{c}$-edge critical non-Hamiltonian graphs when $k \geq 4$. Very recently, Kaemawichanurat and Caccetta~\cite{KC} proved that every $2$-connected $4$-$\gamma_{c}$-edge critical claw-free graph is Hamiltonian. When $k \geq 5$, they found $2$-connected $k$-$\gamma_{c}$-edge critical claw-free graphs which are non-Hamiltonian. For $5 \leq k \leq 6$, they proved that every $k$-$\gamma_{c}$-edge critical claw-free graph is Hamiltonian when it is $3$-connected. So far, the studies on Hamiltonicities of $k$-$\gamma$-vertex critical graphs or $k$-$\gamma_{c}$-vertex critical graphs have not been done.

\indent In this paper, we study Hamiltonicities of $k$-$\gamma$-vertex critical graphs. The results are given in Section 2. We set up notation that are used of in establishing our results in the first subsection. Then, we prove that all $2$-connected $3$-$\gamma$-vertex critical claw-free graphs are Hamiltonian. For $k \geq 3$, we provide a $2$-connected $k$-$\gamma$-vertex critical $K_{1, 4}$-free graph containing a claw as an induced subgraph which is non-Hamiltonian. That is the condition claw-free is necessary to prove that $2$-connected $3$-$\gamma$-vertex critical graphs are Hamiltonian. In Section \ref{HVC}, we present a new method to prove Hamiltonicities of $3$-connected $k$-$\gamma_{c}$-vertex critical claw-free graphs for $3 \leq k \leq 5$. Thus, we set up notation that are used of in establishing our results in the first subsection. We also give $k$-$\gamma_{c}$-vertex critical non-Hamiltonian graphs containing a claw as an induced subgraph for $3 \leq k \leq 5$.

\section{\bf Hamiltonicity of $k$-$\gamma$-vertex critical graphs}\label{HV}
In this section, we prove Hamiltonicities of $k$-$\gamma$-vertex critical graphs. In the first subsection, we introduce a classical method of studying Hamiltonian properties of graphs.

\subsection{\bf Set up and known results}
Suppose $G$ is a $2$-connected non-Hamiltonian graph. Let $C$ be a longest cycle of $G$. We write $\overrightarrow{C}$ to indicate the clockwise orientation of $C$. Similarly, we denote the anticlockwise orientation of $C$ by $\overleftarrow{C}$. In particular, for vertices $u$ and $v$ of $C$ we denote the $(u, v)$-directed segment of $\overrightarrow{C}(\overleftarrow{C})$ by $u\overrightarrow{C}v$ $(u\overleftarrow{C}v)$. The successor (predecessor) of a vertex $v$ of $C$ in $\overrightarrow{C}$ is denoted by $v^{+} (v^{-})$. Furthermore, for $i \geq 1$, $v^{(i + 1)+} = (v^{i+})^{+}$ and $v^{(i + 1)-} = (v^{i-})^{-}$ where $v^{1+} = v^{+}$ and $v^{1-} = v^{-}$. This notation is illustrated by the following figure. Note that we always use an orientation $\overrightarrow{C}$ when we mention about the successor and the predecessor of any vertex of $C$.

\begin{center}
\setlength{\unitlength}{0.8cm}
\begin{picture}(11, 7.5)
\put(5.5, 5){\oval(12, 1.9)}
\put(5.5, 4.8){$\overrightarrow{C}$}

\put(0.5, 4.05){\circle*{0.15}}
\put(1.5, 4.05){\circle*{0.15}}
\put(3.5, 4.05){\circle*{0.15}}
\put(4.5, 4.05){\circle*{0.15}}
\put(5.5, 4.05){\circle*{0.15}}
\put(6.5, 4.05){\circle*{0.15}}
\put(7.5, 4.05){\circle*{0.15}}
\put(9.5, 4.05){\circle*{0.15}}
\put(10.5, 4.05){\circle*{0.15}}

\put(0.5, 5.95){\circle*{0.15}}
\put(1.5, 5.95){\circle*{0.15}}
\put(3.5, 5.95){\circle*{0.15}}
\put(4.5, 5.95){\circle*{0.15}}
\put(5.5, 5.95){\circle*{0.15}}
\put(6.5, 5.95){\circle*{0.15}}
\put(7.5, 5.95){\circle*{0.15}}
\put(9.5, 5.95){\circle*{0.15}}
\put(10.5, 5.95){\circle*{0.15}}

\put(5.5, 6.2){$u$}
\put(4.5, 6.2){$u^{-}$}
\put(3.5, 6.2){$u^{2-}$}
\put(1.5, 6.2){$u^{i-}$}
\put(-0.5, 6.2){$u^{(i + 1)-}$}

\put(6.5, 6.2){$u^{+}$}
\put(7.5, 6.2){$u^{2+}$}
\put(9.5, 6.2){$u^{i+}$}
\put(10.5, 6.2){$u^{(i + 1)+}$}

\multiput(2.5, 6.3)(0.2, 0){3}{\circle*{0.08}}
\multiput(8.5, 6.3)(0.2, 0){3}{\circle*{0.08}}


\put(5.5, 3.5){$v$}
\put(4.5, 3.5){$v^{+}$}
\put(3.5, 3.5){$v^{2+}$}
\put(1.5, 3.5){$v^{i+}$}
\put(-0.5, 3.5){$v^{(i + 1)+}$}

\put(6.5, 3.5){$v^{-}$}
\put(7.5, 3.5){$v^{2-}$}
\put(9.5, 3.5){$v^{i-}$}
\put(10.5, 3.5){$v^{(i + 1)-}$}

\multiput(2.5, 3.7)(0.2, 0){3}{\circle*{0.08}}
\multiput(8.5, 3.7)(0.2, 0){3}{\circle*{0.08}}

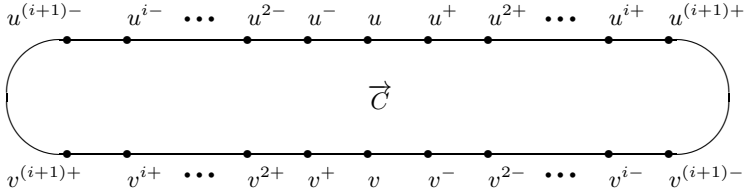
\put(1, 2.5){\small\textbf{Figure 1 :} The $i$th successor and predecessor of $v$.}
\end{picture}
\end{center}
\vskip -51 pt

\noindent Let $H$ be a component of $G - C$ and $X = N_{C}(H)$. Suppose $|X| = d$. We may order the vertices of the set $X$ as $x_{1}, x_{2}, ..., x_{d}$ according to $\overrightarrow{C}$. We, further, let
\vskip 5 pt

\indent $X^{+} = \{x^{+}_{1}, x^{+}_{2}, ..., x^{+}_{d}\}$ and $X^{-} = \{x^{-}_{1}, x^{-}_{2}, ..., x^{-}_{d}\}$.
\vskip 5 pt

\noindent For any vertices $u$ and $v$ of $C$, if $|u\overrightarrow{C}v| = t$, then we can order the vertices in $u\overrightarrow{C}v$ as $u, u^{+}, u^{2+}, ...,$ $u^{(t - 1)+}$ such that $u^{(t - 1)+} = v$. We let $\overrightarrow{C}[u, v] = \{u, u^{+}, u^{2+}, ..., u^{(t - 1)+}\}$. All subscripts are taken modulo $d$ throughout. Favaron et al.~\cite{FTZ} provided some structural properties described in Lemmas \ref{lem 21} and \ref{lem 22}. Note that, for any two vertices $u, v$ and subgraph $H$ of $G$, $uP_{H}v$ denotes a path from $u$ to $v$ which the internal vertices are in $V(H)$.
\vskip 5 pt

\begin{lem}\label{lem 21}\cite{FTZ}
$X^{+} \cap X = \emptyset$ and $X^{-} \cap X = \emptyset$.
\end{lem}

\begin{lem}\label{lem 22}\cite{FTZ}
For any two vertices $u, v \in X^{+}$ (or $u, v \in X^{-}$), then there is no $uP_{G - C}v$ path, in particular, $uv \notin E(G)$.
\end{lem}

\noindent Lemmas \ref{lem h0} and \ref{lem h0n} are well known (see Brousek~\cite{B2}) and proved under the condition that $G$ is claw-free non-Hamiltonian graph.
\vskip 5 pt

\begin{lem}\label{lem h0}\cite{B2}
For all $1 \leq i \leq d$, $x^{+}_{i}x^{-}_{i} \in E(G)$.
\end{lem}

\begin{lem}\label{lem h1}\cite{B2}
For $1 \leq i \neq j \leq d$, $\{x_{i}x^{+}_{j}, x_{i}x^{2+}_{j}, x^{+}_{i}x^{2+}_{j}, x_{i}x^{-}_{j}, x_{i}x^{2-}_{j}, x^{-}_{i}x^{2-}_{j}\} \cap E(G) = \emptyset$.
\end{lem}

\begin{lem}\label{lem h0n}\cite{B2}
For all $1 \leq i \leq d$, $|\overrightarrow{C}[x^{+}_{i}, x^{-}_{i + 1}]| \geq 3$.
\end{lem}
\vskip 5 pt

\indent The study $k$-$\gamma$-vertex critical graphs has been started by Brigham et al.~\cite{BCD} and continued by Fulman et al.~\cite{FHM}. Moreover, Ananchuen and Plummer~\cite{AP} established some observations on these graphs.

\begin{lem}\label{lem 1}\cite{AP}
For $k \geq 2$, let $G$ be a $k$-$\gamma$-vertex critical graph. Moreover, for a vertex $v$ of $G$, we let $D_{v}$ be a $\gamma$-set of $G - v$. Then
\vskip 5 pt

\indent (i) $|D_{v}| = k - 1$ and
\vskip 5 pt

\indent (ii) every vertex in $D_{v}$ is not adjacent to $v$.
\end{lem}
\vskip 5 pt

\subsection{\bf Main results}
In this subsection, we prove that all $3$-$\gamma$-vertex critical claw free graphs are Hamiltonian. We first give a construction of $k$-$\gamma$-vertex critical graphs $G_{k}$ which are non-Hamiltonian for $k \geq 3$. Let $C_{k} = b_{1}b_{2}...b_{k}b_{1}$ be a cycle of length $k$. We construct the graph $S_{k}$ by subdividing the edge $b_{i}b_{i + 1}$ with a vertex $a_{i + 1}$ for every $1 \leq i \leq k - 1$, subdividing the edge $b_{k}b_{1}$ with a vertex $a_{1}$ and then, joining $a_{i}$ to $a_{i + 1}$ for $1 \leq i \leq k - 1$ and joining $a_{k}$ to $a_{1}$. Our graph $G_{k}$ is constructed from $S_{k}$ by subdividing the edge $a_{i}a_{i + 1}$ with a vertex $c_{i}$ for all $1 \leq i \leq k - 1$, subdividing the edge $a_{k}a_{1}$ with a vertex $c_{k}$ and finally, joining $c_{i}$ to $c_{i + 1}$ for all $1 \leq i \leq k - 1$ and joining $c_{k}$ to $c_{1}$. The following figures illustrated $G_{3}$ and $G_{8}$ for examples. It is worth noting that, when $k = 3$, the graph $G_{3}$ has been found earlier by Ananchuen and Plummer~\cite{AP}.
\vskip 20 pt

\setlength{\unitlength}{1cm}
\begin{center}
\begin{picture}(6,6)
\put(0, 0){\circle*{0.15}}
\put(-0.5, -0.2){\small{$b_{3}$}}

\put(6, 0){\circle*{0.15}}
\put(6.2, -0.2){\small{$b_{2}$}}

\put(3, 6){\circle*{0.15}}
\put(3, 6.2){\small{$b_{1}$}}

\put(3, 0){\circle*{0.15}}
\put(3, -0.3){\small{$a_{3}$}}

\put(2.25, 1.5){\circle*{0.15}}
\put(1.8, 1.5){\small{$c_{3}$}}

\put(3.75, 1.5){\circle*{0.15}}
\put(3.95, 1.5){\small{$c_{2}$}}

\put(1.5, 3){\circle*{0.15}}
\put(1.2, 3.2){\small{$a_{1}$}}

\put(3, 3){\circle*{0.15}}
\put(3, 3.2){\small{$c_{1}$}}

\put(4.5, 3){\circle*{0.15}}
\put(4.6, 3.2){\small{$a_{2}$}}

\put(0, 0){\line(1, 0){6}}
\put(0, 0){\line(1, 2){3}}
\put(6, 0){\line(-1, 2){3}}

\put(1.5, 3){\line(1, 0){3}}
\put(1.5, 3){\line(1, -2){1.5}}
\put(4.5, 3){\line(-1, -2){1.5}}

\put(3, 3){\line(1, -2){0.75}}
\put(3, 3){\line(-1, -2){0.75}}
\put(2.25, 1.5){\line(1, 0){1.5}}

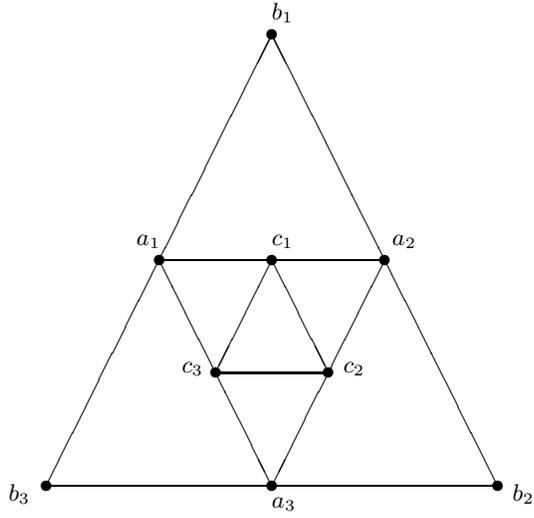
\put(1.5, -1.5){\textbf{Figure 2 :} The graph $G_{3}$}
\end{picture}
\end{center}
\vskip 10 pt

\setlength{\unitlength}{1cm}
\begin{center}
\begin{picture}(5,10)
\put(1, 0){\circle*{0.15}}
\put(4, 0){\circle*{0.15}}
\put(-1.2, 2.18){\circle*{0.15}}
\put(-1.2, 5.18){\circle*{0.15}}
\put(1, 7.38){\circle*{0.15}}
\put(0.7, 7.52){\small{$b_{8}$}}

\put(4, 7.38){\circle*{0.15}}
\put(4, 7.52){\small{$b_{1}$}}

\put(6.2, 2.18){\circle*{0.15}}
\put(6.2, 5.18){\circle*{0.15}}
\put(6.4, 5.18){\small{$b_{2}$}}
\put(6.4, 4.18){\small{$\vdots$}}

\put(2.5, 0){\circle*{0.15}}
\put(0, 1){\circle*{0.15}}
\put(-1.2, 3.4){\circle*{0.15}}
\put(-0.05, 6.35){\circle*{0.15}}
\put(-0.5, 6.4){\small{$a_{8}$}}

\put(2.5, 7.36){\circle*{0.15}}
\put(2.4, 7.52){\small{$a_{1}$}}

\put(5, 6.35){\circle*{0.15}}
\put(5.2, 6.4){\small{$a_{2}$}}

\put(6.2, 3.9){\circle*{0.15}}
\put(5, 1){\circle*{0.15}}

\put(1.2, 0.5){\circle*{0.15}}
\put(-0.5, 1.96){\circle*{0.15}}
\put(-0.55, 5.1){\circle*{0.15}}
\put(1.16, 6.81){\circle*{0.15}}
\put(1.16, 6.4){\small{$c_{8}$}}

\put(3.86, 6.81){\circle*{0.15}}
\put(3.58, 6.4){\small{$c_{1}$}}

\put(5.6, 5.1){\circle*{0.15}}
\put(5, 5.18){\small{$c_{2}$}}
\put(5, 4.18){\small{$\vdots$}}

\put(5.6, 2.5){\circle*{0.15}}
\put(3.65, 0.5){\circle*{0.15}}

\put(1, 0){\line(1, 0){3}}
\put(1, 0){\line(-1, 1){2.18}}
\put(4, 0){\line(1, 1){2.18}}
\put(-1.2, 2.18){\line(0, 1){3}}
\put(6.2, 2.18){\line(0, 1){3}}
\put(-1.2, 5.18){\line(1, 1){2.18}}
\put(6.2, 5.18){\line(-1, 1){2.18}}
\put(1, 7.38){\line(1, 0){3}}

\put(2.5, 0){\line(5, 2){2.5}}
\put(2.5, 0){\line(-5, 2){2.5}}
\put(0, 1){\line(-1, 2){1.2}}
\put(-1.2, 3.45){\line(2, 5){1.15}}
\put(-0.05, 6.35){\line(5, 2){2.57}}
\put(2.45, 7.38){\line(5, -2){2.57}}
\put(5, 6.3){\line(1, -2){1.2}}
\put(6.2, 4){\line(-2, -5){1.2}}

\put(1.25, 0.5){\line(1, 0){2.4}}
\put(1.25, 0.5){\line(-6, 5){1.7}}
\put(-0.5, 1.9){\line(0, 1){3.2}}
\put(-0.5, 5.1){\line(1, 1){1.62}}
\put(1.23, 6.8){\line(1, 0){2.6}}
\put(3.85, 6.8){\line(1, -1){1.7}}
\put(5.6, 5.1){\line(0, -1){2.5}}
\put(5.65, 2.55){\line(-1, -1){2}}

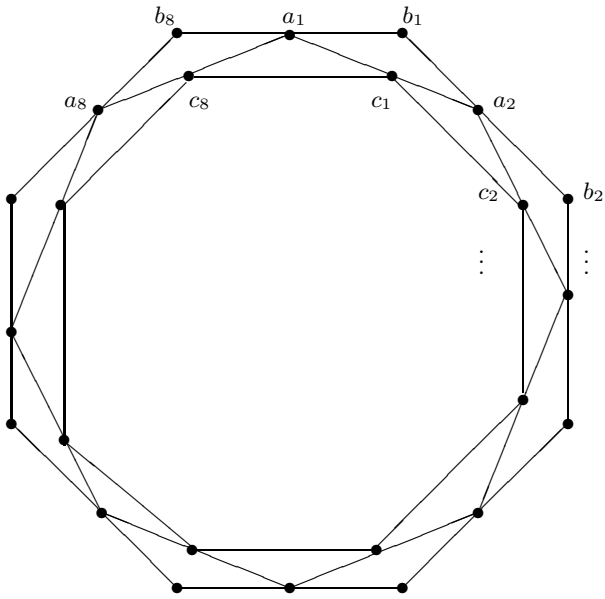
\put(1, -1.3){\textbf{Figure 3 :} The graph $G_{8}$}
\end{picture}
\end{center}
\vskip 40 pt

\noindent  In the following, we let
\vskip 5 pt

\indent $A_{i} = \{a_{i}, b_{i}, c_{i}\}$ for $1 \leq i \leq k$,
\vskip 5 pt

\indent $A = \{a_{1}, a_{2}, ..., a_{k}\}$ and $B = \{b_{1}, b_{2}, ..., b_{k}\}$.
\vskip 5 pt

\noindent The following lemma gives that $\gamma(G_{k}) = k$.
\vskip 5 pt

\begin{lem}\label{lem gk1}
For a positive integer $k \geq 3$, we have that $\gamma(G_{k}) = k$.
\end{lem}

\proof
Clearly, $A \succ G_{k}$. By the minimality of $\gamma(G_{k})$, we have $\gamma(G_{k}) \leq k$.

We will show that $\gamma(G_{k}) \geq k$. Suppose that $\gamma(G) < k$ and let $D$ be a $\gamma$-set such that $|D \cap A|$ is maximum.

If $D \cap A = \emptyset$, then $B \subseteq D$ to dominate $B$. Thus $\gamma(G_{k}) = |D| \geq |B| = k$, a contradiction. Hence, $D \cap A \neq \emptyset$. Since $|D| < k$, it follows that $A \nsubseteq D$. Thus, there exists $i$ such that $a_{i - 1} \in D$ but $a_{i} \notin D$. To dominate $a_{i}$, we have that $\{b_{i - 1}, c_{i - 1}, b_{i}, c_{i}\} \cap D \neq \emptyset$. If $b_{i - 1} \in D$, then $D' = (D - \{b_{i - 1}\}) \cup \{a_{i}\} \succ G_{k}$ but $|D' \cap A| > |D \cap A|$, contradicting the maximality of $|D \cap A|$. Similarly, if $c_{i - 1} \in D$, then $D' = (D - \{c_{i - 1}\}) \cup \{a_{i}\} \succ G_{k}$ but $|D' \cap A| > |D \cap A|$, a contradiction. Hence, $b_{i} \in D$ or $c_{i} \in D$.

We first assume that $c_{i} \in D$. To dominate $b_{i}$, we have $b_{i} \in D$ or $a_{i + 1} \in D$. If $a_{i + 1} \in D$, then $D' = (D - \{c_{i}\}) \cup \{a_{i}\} \succ G_{k}$ but $|D' \cap A| > |D \cap A|$, a contradiction. Hence, $b_{i} \in D$. In this case, $D' = (D - \{b_{i}\}) \cup \{a_{i}\} \succ G_{k}$ but $|D' \cap A| > |D \cap A|$ which is a contradiction.

Thus, we may assume that $c_{i} \notin D$. Therefore, $b_{i} \in D$. To dominate $c_{i}$, we have $a_{i + 1} \in D$ or $c_{i + 1} \in D$. In both cases, $D' = (D - \{b_{i}\}) \cup \{a_{i}\} \succ G_{k}$ but $|D' \cap A| > |D \cap A|$, a contradiction.

Therefore, $\gamma(G_{k}) \geq k$ implying that $\gamma(G_{k}) = k$. This completes the proof.
\qed

\vskip 5 pt

\indent The next lemma shows that, for $k \geq 3$, the graph $G_{k}$ is $k$-$\gamma$-vertex critical which is non-Hamiltonian.
\vskip 5 pt

\begin{lem}\label{lem A}
For $k \geq 3$, the graph $G_{k}$ is $k$-$\gamma$-vertex critical $K_{1, 4}$-free containing $K_{1, 3}$ which is non-Hamiltonian.
\end{lem}
\proof
In view of Lemma \ref{lem gk1}, we have $\gamma(G_{k}) = k$. We will establish the criticality. For $1 \leq i \leq k$, we have
\vskip 5 pt

\indent $D_{a_{i}} = A - \{a_{i}\}, D_{b_{i}} = (A - \{a_{i}, a_{i + 1}\})\cup \{c_{i}\}$ and $D_{c_{i}} = (A - \{a_{i}, a_{i + 1}\})\cup \{b_{i}\}$
\vskip 5 pt

\noindent which the subscripts are taken modulo $k$. We see that $|D_{v}| = k - 1$ for all $v \in V(G_{k})$. Hence, $G_{k}$ is a $k$-$\gamma$-vertex critical. It is not difficult to see that $G_{k}$ is a $K_{1, 4}$-free graph containing a claw $G_{k}[\{a_{1}, b_{k}, b_{1}, c_{1}\}]$ centered at $a_{1}$ as an induced subgraph.
\vskip 5 pt

\indent Finally, suppose that the graph $G_{k}$ is Hamiltonian. We observe that the edges incident with the degree $2$ vertices, $b_{1}, b_{2}, ..., b_{k}$ must be in a Hamiltonian cycle. These edges induce a cycle which does not contain the vertices $c_{1}, c_{2}, ..., c_{k}$, a contradiction. Hence, $G_{k}$ is non-Hamiltonian and this completes the proof.
\qed


\vskip 5 pt

\indent We will show that every $3$-$\gamma$-vertex critical graph is Hamiltonian when it is claw-free.
\vskip 5 pt

\begin{thm}\label{thm A}
Let $G$ be a $2$-connected $3$-$\gamma$-vertex critical claw-free graph. Then $G$ is Hamiltonian.
\end{thm}
\proof
Suppose to the contrary that $G$ does not contain a Hamiltonian cycle. Let $C$ be a longest cycle of $G$. Thus there exists a component $H$ of $G - C$. We let, further, that $X = N_{C}(H)$. We may order the vertices in $X$ as $x_{1}, ..., x_{d}$ according to $\overrightarrow{C}$. Because $G$ is $2$-connected, it follows that $|X| \geq 2$. We note by Lemma \ref{lem h0n} that $|\overrightarrow{C}[x^{+}_{i}, x^{-}_{i + 1}]| \geq 3$ for all $1 \leq i \leq d$. We need to establish the following claims.
\vskip 5 pt

\noindent \textbf{Claim 1 :} For a vertex $y \in V(C)$, we have that $D_{y} \cap V(C) \neq \emptyset$.
\vskip 5 pt

\indent Let $y \in V(C)$. Suppose to the contrary that $D_{y} \cap V(C) = \emptyset$. To dominate $H$, $D_{y} \cap V(H) \neq \emptyset$. Let $u \in D_{y} \cap V(H)$. Thus, $|D_{y} - \{u\}| = 1$ by Lemma \ref{lem 1}($i$). Let $\{v\} = D_{y} - \{u\}$. Lemma \ref{lem 21} yields that $u$ is not adjacent to any vertex in $X^{+} \cup X^{-}$. Since $D_{y} \succ G - y$ , it follows that $vx^{+}_{1}, vx^{+}_{2} \in E(G)$ when $y \in \{x^{-}_{1}, x^{-}_{2}\}$ or $vx^{-}_{1}, vx^{-}_{2} \in E(G)$ otherwise. This contradicts Lemma \ref{lem 22}. Therefore $D_{y} \cap V(C) \neq \emptyset$ and thus, establishing Claim 1.
\vskip 5 pt

\noindent \textbf{Claim 2 :} For a vertex $y \in V(C)$, we have that $|D_{y} \cap X| \leq 1$.
\vskip 5 pt

\indent Suppose to the contrary that $|D_{y} \cap X| \geq 2$. By Lemma \ref{lem 1}($i$), $|D_{y} \cap X| = 2$. Let $D_{y}  = \{x_{i}, x_{j}\}$. Without loss of generality, let $y \notin \overrightarrow{C}[x_{i}, x_{j}]$. Thus $\{x_{i}, x_{j}\} \succ \overrightarrow{C}[x_{i}, x_{j}]$. We note by Lemma \ref{lem h1} that $x_{i}$ does not dominate $\overrightarrow{C}[x_{i}, x^{2-}_{j}]$. Let $|\overrightarrow{C}[x_{i}, x^{2-}_{j}]| = t_{0}$. We can order the vertices in $\overrightarrow{C}[x_{i}, x^{2-}_{j}]$ to be $x_{i}, x^{+}_{i}, ..., x^{(t_{0} - 1)+}_{i}$ where $x^{(t_{0} - 1)+}_{i} = x^{-2}_{j}$. Then, we let $r_{0} = max\{r : 1 \leq r < t_{0}$ and $x_{i}x^{r+}_{i} \in E(G)\}$. Thus $x_{j}x^{(r_{0} + 1)+}_{i} \in E(G)$. We see that
\vskip 5 pt

\indent $x_{i}x^{r_{0}}_{i}\overleftarrow{C}x^{+}_{i}x^{-}_{i}\overleftarrow{C}x^{+}_{j}x^{-}_{j}\overleftarrow{C}x^{(r_{0} + 1)+}_{i}x_{j}P_{H}x_{i}$
\vskip 5 pt

\noindent is a cycle longer than $C$, a contradiction. Hence, $|D_{y} \cap X| \leq 1$ and we establish Claim 2.
\vskip 12 pt

\indent Consider $G - x_{1}$. In view of Lemma \ref{lem 1}($i$), $|D_{x_{1}}| = 2$. We distinguish two cases.
\vskip 5 pt

\noindent \textbf{Case 1 :} $D_{x_{1}} \cap X = \emptyset$.
\vskip 5 pt

\indent This implies by Claim 1 that $D_{x_{1}} \cap V(C) \neq \emptyset$. To dominate $H$, $D_{x_{1}} \cap V(H) \neq \emptyset$. Thus, $|D_{x_{1}} \cap (V(C) - X)| = 1$ and $|D_{x_{1}} \cap V(H)| = 1$. Let $\{a\} = D_{x_{1}} \cap (V(C) - X)$ and $\{b\} = D_{x_{1}} \cap V(H)$. Clearly, $b$ is not adjacent to any vertex in $V(G) - (X \cup (V(H)))$. Therefore, $a \succ V(G) - (X \cup (V(H)))$.
\vskip 5 pt

\indent Consider $G - a$. We have by Claim 1 that $D_{a} \cap V(C) \neq \emptyset$. Since $a \succ V(G) - (X \cup (V(H)))$, by Lemma \ref{lem 1}($ii$), $D_{a} \cap (V(C) - X) = \emptyset$. Thus $D_{a} \subseteq X \cup V(H)$. If $D_{a} \subseteq V(H)$, then, by Lemma \ref{lem 21}, $D_{a}$ does not dominate $X^{+}$, a contradiction. Thus, $D_{a} \cap X \neq \emptyset$. By Claim 2, $|D_{a} \cap X| = 1$, moreover, $|D_{a} \cap V(H)| = 1$. Let $\{x_{i}\} = D_{a} \cap X$ and $\{h\} = D_{a} \cap V(H)$. By Lemma \ref{lem 21}, $h$ is not adjacent to any vertex in $V(C) - X$. This implies that $x_{i} \succ V(C) - X$. In particular, $x_{i}x^{+}_{j} \in E(G)$ for some $x_{j} \in X - \{x_{i}\}$. This contradicts Lemma \ref{lem h1}. Hence, Case 1 cannot occur.
\vskip 5 pt

\noindent \textbf{Case 2 :} $D_{x_{1}} \cap X \neq \emptyset$.
\vskip 5 pt

\indent This implies by Claim 2 that $|D_{x_{1}} \cap X| = 1$. Thus $|D_{x_{1}} - X| = 1$. Let $\{x_{i}\} = D_{x_{1}} \cap X$ and $\{u\} = D_{x_{1}} - X$. In view of Lemma \ref{lem h1}, $x_{i}$ is not adjacent to any vertex in $W = \{x^{+}_{1}, x^{2+}_{1}, x^{-}_{1}, x^{2-}_{1}\}$. Thus $u \succ W$. This implies by Lemma \ref{lem 21} that $u \notin V(H)$. In the following,  we assume by Lemma \ref{lem h0} that $x^{+}_{i}x^{-}_{i}, x^{+}_{1}x^{-}_{1} \in E(G)$.
\vskip 5 pt

\indent Suppose that $u \notin V(C)$. Let $|\overrightarrow{C}[x^{+}_{i}, x^{-}_{1}]| = t_{1}$. The vertices in $\overrightarrow{C}[x^{+}_{i}, x^{-}_{1}]$ can be ordered as $x^{+}_{i}, x^{2+}_{i}, ..., x^{t_{1}+}_{i}$ where $x^{t_{1}+}_{i} = x^{-}_{1}$. Then we let $r_{1} = max\{r : 1 \leq r < t_{1}$ and $x_{i}x^{r+}_{i} \in E(G)\}$. Therefore, $ux^{(r_{1} + 1)+}_{i}$. Clearly,
\vskip 5 pt

\indent $x_{i}x^{r_{1}+}_{i}\overleftarrow{C}x^{+}_{i}x^{-}_{i}\overleftarrow{C}x^{+}_{1}ux^{(r_{1} + 1)+}_{i}\overrightarrow{C}x_{1}P_{H}x_{i}$
\vskip 5 pt

\noindent is a cycle longer than $C$, a contradiction. Therefore, $u \in V(C)$. Thus either $u \in \overrightarrow{C}[x^{+}_{i}, x^{-}_{1}]$ or $u \in \overrightarrow{C}[x^{+}_{1}, x^{-}_{i}]$.
\vskip 5 pt

\noindent \textbf{Subcase 2.1 :} $u \in \overrightarrow{C}[x^{+}_{i}, x^{-}_{1}]$.
\vskip 5 pt

\indent Recall that $ux^{+}_{1}, ux^{2+}_{1} \in E(G)$ but $x_{i}x^{+}_{1}, x_{i}x^{2+}_{1} \notin E(G)$. Thus, Lemma \ref{lem 22} implies that $u  \neq x^{+}_{i}$. Let $|\overrightarrow{C}[x^{+}_{1}, x^{-}_{i}]| = t_{2}$. We can order the vertices in $\overrightarrow{C}[x^{+}_{1}, x^{-}_{i}]$ as $x^{+}_{1}, x^{2+}_{1}, ..., x^{t_{2}+}_{1}$ where $x^{t_{2}+}_{1} = x^{-}_{i}$. Then we let $r_{2} = min\{r : 1 < r \leq t_{2}$ and $x_{i}x^{r+}_{1} \in E(G)\}$. By Lemma \ref{lem h1}, $r_{2} > 2$. We then have that $ux^{(r_{2} - 1)+}_{1} \in E(G)$ because $\{x_{i}, u\} \succ G - x_{1}$. If $u = x^{-}_{1}$, then
\vskip 5 pt

\indent $x_{i}x^{r_{2}+}_{1}\overrightarrow{C}x^{-}_{i}x^{+}_{i}\overrightarrow{C}ux^{(r_{2} - 1)+}_{1}\overleftarrow{C}x_{1}P_{H}x_{i}$
\vskip 5 pt

\noindent is a cycle longer than $C$, a contradiction. Hence, suppose that $u \neq x^{-}_{1}$. We distinguish two subcases.
\vskip 5 pt

\noindent \textbf{Subcase 2.1.1 :} $u^{+}u^{-} \in E(G)$.
\vskip 5 pt

\indent Since $u \succ W$, $ux^{-}_{1} \in E(G)$. Thus,
\vskip 5 pt

\indent $x_{i}x^{r_{2}+}_{1}\overrightarrow{C}x^{-}_{i}x^{+}_{i}\overrightarrow{C}u^{-}u^{+}\overrightarrow{C}x^{-}_{1}ux^{(r_{2} - 1)+}_{1}\overleftarrow{C}x_{1}P_{H}x_{i}$
\vskip 5 pt

\noindent is  cycle longer than $C$, a contradiction. Thus this subcase cannot occur.
\vskip 5 pt

\noindent \textbf{Subcase 2.1.2 :} $u^{+}u^{-} \notin E(G)$
\vskip 5 pt

\indent By claw-freeness of $G$, $x^{(r_{2} - 1)+}_{1}$ is adjacent to $u^{-}$ or $u^{+}$. If $x^{(r_{2} - 1)+}_{1}u^{-} \in E(G)$, then
\vskip 5 pt

\indent $x_{i}x^{r_{2}+}_{1}\overrightarrow{C}x^{-}_{i}x^{+}_{i}\overrightarrow{C}u^{-}x^{(r_{2} - 1)+}_{1}\overleftarrow{C}x^{+}_{1}u\overrightarrow{C}x_{1}P_{H}x_{i}$
\vskip 5 pt

\noindent is a cycle longer than $C$, a contradiction. Thus $x^{(r_{2} - 1)+}_{1}u^{+} \in E(G)$. We see that
\vskip 5 pt

\indent $x_{i}x^{r_{2}+}_{1}\overrightarrow{C}x^{-}_{i}x^{+}_{i}\overrightarrow{C}ux^{+}_{1}\overrightarrow{C}x^{(r_{2} - 1)+}_{1}u^{+}\overrightarrow{C}x_{1}P_{H}x_{i}$
\vskip 5 pt

\noindent is a cycle longer than $C$, a contradiction. Therefore, Subcase 2.1 cannot occur.
\vskip 10 pt

\noindent \textbf{Subcase 2.2 :} $u \in \overrightarrow{C}[x^{+}_{1}, x^{-}_{i}]$.
\vskip 5 pt

\indent Since $u \succ W$, we have that $ux^{-}_{1}, ux^{2-}_{1} \in E(G)$. By Lemma \ref{lem 22}, $u \neq x^{-}_{i}$. Let $|\overrightarrow{C}[x^{+}_{i}, x^{-}_{1}]| = t_{3}$. The vertices in $\overrightarrow{C}[x^{+}_{i}, x^{-}_{1}]$ can be ordered as $x^{+}_{i}, x^{2+}_{i}, ..., x^{t_{3}+}_{i}$ where $x^{t_{3}+}_{i} = x^{-}_{1}$. Then we let $r_{3} = max\{r : 1 \leq r < t_{3}$ and $x_{i}x^{r+}_{i} \in E(G)\}$. Because $\{u, x_{i}\} \succ G - x_{1}$ and $x_{1} \notin \overrightarrow{C}[x^{+}_{i}, x^{-}_{1}]$, it follows that $ux^{(r_{3} + 1)+}_{i}$. If $u = x^{+}_{1}$, then
\vskip 5 pt

\indent $x_{i}x^{r_{3}+}_{i}\overleftarrow{C}x^{+}_{i}x^{-}_{i}\overleftarrow{C}ux^{(r_{3} + 1)+}_{i}\overrightarrow{C}x_{1}P_{H}x_{i}$
\vskip 5 pt

\noindent is a cycle longer than $C$, a contradiction. Thu, $u \neq x^{+}_{1}$. We distinguish two more subcases.
\vskip5 pt

\noindent \textbf{Subcase 2.2.1 :} $u^{+}u^{-} \in E(G)$.
\vskip 5 pt

\indent As $u \succ W$, we must have $ux^{+}_{1} \in E(G)$. Thus,
\vskip 5 pt

\indent $x_{i}x^{r_{3}+}_{1}\overleftarrow{C}x^{+}_{i}x^{-}_{i}\overleftarrow{C}u^{+}u^{-}\overleftarrow{C}x^{+}_{1}ux^{(r_{3} + 1)+}_{i}\overrightarrow{C}x_{1}P_{H}x_{i}$
\vskip 5 pt

\noindent is  cycle longer than $C$, a contradiction. Thus this subcase cannot occur.
\vskip 5 pt

\noindent \textbf{Subcase 2.2.2 :} $u^{+}u^{-} \notin E(G)$
\vskip 5 pt

\indent By claw-freeness of $G$, $x^{(r_{3} + 1)+}_{i}$ is adjacent to $u^{-}$ or $u^{+}$. If $x^{(r_{3} + 1)+}_{i}u^{-} \in E(G)$, then
\vskip 5 pt

\indent $x_{i}x^{r_{3}+}_{i}\overleftarrow{C}x^{+}_{i}x^{-}_{i}\overleftarrow{C}ux^{-}_{1}\overleftarrow{C}x^{(r_{3} + 1)+}_{i}u^{-}\overleftarrow{C}x_{1}P_{H}x_{i}$
\vskip 5 pt

\noindent is a cycle longer than $C$, a contradiction. Thus $x^{(r_{3} + 1)+}_{i}u^{+} \in E(G)$. We see that
\vskip 5 pt

\indent $x_{i}x^{r_{3}+}_{i}\overleftarrow{C}x^{+}_{i}x^{-}_{i}\overleftarrow{C}u^{+}x^{(r_{3} + 1)+}_{i}\overrightarrow{C}x^{-}_{1}u\overleftarrow{C}x_{1}P_{H}x_{i}$
\vskip 5 pt

\noindent is a cycle longer than $C$, a contradiction. Therefore, Subcase 2.2 cannot occur. Hence, $G$ is Hamiltonian and this completes the proof.
\qed
\vskip 5 pt

\indent In view of Lemma \ref{lem A}, we see that the condition claw-free in Theorem \ref{thm A} is best possible.
\vskip 5 pt

\section{\bf Hamiltonicity of $k$-$\gamma_{c}$-vertex critical graphs}\label{HVC}
In this section we present a new method (without using the classical method) to prove that every $2$-connected $3$-$\gamma_{c}$-vertex critical claw-free graphs is Hamiltonian. Moreover, for $4 \leq k \leq 5$, we prove that every $3$-connected $k$-$\gamma_{c}$-vertex critical claw-free graph is Hamiltonian. For $3 \leq k \leq 5$, we find some $2$-connected $k$-$\gamma_{c}$-vertex critical non-Hamiltonian graphs containing a claw as an induced subgraph. We first give some related results that we use in the first subsection.

\subsection{\bf Set up and known results}\label{SS2}
In this subsection, we begin with a result of Chv\'{a}tal~\cite{C} which is a property of a graph when it is Hamiltonian.
\vskip 5 pt

\begin{lem}\label{lem Ch}
Let $G$ be a Hamiltonian graph and $S \subseteq V(G)$ a cut set of $G$. Then $\frac{|S|}{\omega(G - S)} \geq 1$.
\end{lem}

\indent In the study on Hamiltonicity of claw-free graphs, there is a technique the so called \emph{local completion} introduced by Ryj\'{a}$\check{c}$ek~\cite{R}. The purpose of this technique is to find the closure $cl(G)$ of a claw-free graph $G$. As we use only the results of this operation, we may omit the detail of finding the closure of a graph. However, it is worth noting that $V(G) = V(cl(G))$ and $E(G) \subseteq E(cl(G))$. Brousek et al.~\cite{BRF} used this closure operation to establish the Hamiltonian properties of $\{K_{1, 3}, N_{s_{1}, s_{2}, s_{3}}\}$-free graphs. The following graphs were found in \cite{BRF}.
\vskip 10 pt

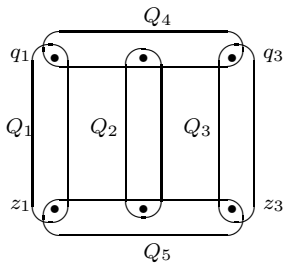
\noindent \textbf{The Class} $\mathcal{F}_{1} :$\\
\indent Let $Q_{1}, ..., Q_{5}$ be complete graphs of order at least three. For $1 \leq i \leq 3$, let $q_{i}, z_{i}$ be two different vertices of $Q_{i}$. Moreover, let $q'_{1}, q'_{2}, q'_{3}$ be three different vertices of $Q_{4}$ and $z'_{1}, z'_{2}, z'_{3}$ be three different vertices of $Q_{5}$. A graph in this class is constructed from $Q_{1}, ..., Q_{5}$ by identifying $q'_{i}$ with $q_{i}$ and $z'_{i}$ with $z_{i}$ for $1 \leq i \leq 3$. A graph in this class is illustrated by Figure 4(a). Note that, an \emph{oval} in the figure denotes a compete graph.
\vskip 10 pt

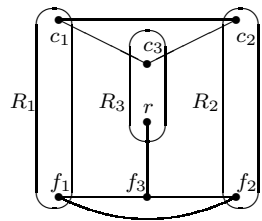
\noindent \textbf{The Class} $\mathcal{F}_{2} :$\\
\indent Let $c_{1}c_{2}c_{3}c_{1}$ and $f_{1}f_{2}f_{3}f_{1}$ be two disjoint triangles. We also let $R_{1}$ and $R_{2}$ be two complete graphs of order at least three and $R_{3}$ a complete graph of order at least two Let $c'_{i}, f'_{i}$ be two different vertices of $R_{i}$ for $1 \leq i \leq 2$ and let $c'_{3}, r$ be two different vertices of $R_{3}$. A graph in this class is obtained by identifying $c'_{i}$ with $c_{i}$ and $f'_{i}$ with $f_{i}$ for $1 \leq i \leq 2$ and identifying $c'_{3}$ with $c_{3}$ and adding an edge $rf_{3}$. A graph in this class is illustrated by Figure 4(b).
\vskip 10 pt

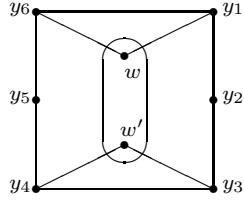
\noindent \textbf{The Class} $\mathcal{F}_{3} :$\\
\indent Let $y_{1}y_{2}...y_{6}y_{1}$ be a cycle of six vertices and $K$ a complete graph of order at least three. Let $w$ and $w'$ be two different vertices of $K$. We define a graph $G$ in the class $\mathcal{F}_{3}$ by adding edges $wy_{1}, wy_{6}, w'y_{3}, w'y_{4}$. A graph in this class is illustrated by Figure 4(c).
\vskip 2 pt

\setlength{\unitlength}{0.59cm}
\begin{center}
\begin{picture}(10,5.9)

\put(-1, 0.3){\circle*{0.15}}
\put(-2, 0.3){\footnotesize{$z_{1}$}}
\put(-1, 3.7){\circle*{0.15}}
\put(-2, 3.7){\footnotesize{$q_{1}$}}
\put(-1.1, 2){\oval(0.8, 4)}
\put(-2.1, 2){\footnotesize{$Q_{1}$}}

\put(1, 0.3){\circle*{0.15}}
\put(1, 3.7){\circle*{0.15}}
\put(1, 2){\oval(0.8, 3.8)}
\put(-0.2, 2){\footnotesize{$Q_{2}$}}

\put(3, 0.3){\circle*{0.15}}
\put(3.7, 0.3){\footnotesize{$z_{3}$}}
\put(3, 3.7){\circle*{0.15}}
\put(3.7, 3.7){\footnotesize{$q_{3}$}}
\put(3.1, 2){\oval(0.8, 4)}
\put(1.9, 2){\footnotesize{$Q_{3}$}}

\put(1, 3.9){\oval(4.5, 0.8)}
\put(1, 0.1){\oval(4.5, 0.8)}
\put(1, 4.5){\footnotesize{$Q_{4}$}}
\put(1, -0.8){\footnotesize{$Q_{5}$}}

\put(-3, -2){\small{\textbf{Figure 4(a) : } The Class $\mathcal{F}_{1}$}}
\put(6, -2){\small{\textbf{Figure 4(b) :} The Class $\mathcal{F}_{2}$}}

\put(7, 0){\circle*{0.15}}
\put(7, 4){\circle*{0.15}}
\put(11, 0){\circle*{0.15}}
\put(11, 4){\circle*{0.15}}
\put(9, 3){\circle*{0.15}}
\put(9, 1.7){\circle*{0.15}}
\put(9, 0){\circle*{0.15}}

\put(5.9, 2){\footnotesize{$R_{1}$}}
\put(7.9, 2){\footnotesize{$R_{3}$}}
\put(10, 2){\footnotesize{$R_{2}$}}
\put(8.9, 1.9){\footnotesize{$r$}}
\put(8.5, 0.2){\footnotesize{$f_{3}$}}

\put(6.8, 3.5){\footnotesize{$c_{1}$}}
\put(11, 3.5){\footnotesize{$c_{2}$}}
\put(8.9, 3.3){\footnotesize{$c_{3}$}}
\put(6.8, 0.2){\footnotesize{$f_{1}$}}
\put(11, 0.2){\footnotesize{$f_{2}$}}

\put(7, 0){\line(1, 0){4}}
\put(7, 4){\line(2, -1){2}}
\put(11, 4){\line(-2, -1){2}}
\put(11, 4){\line(-1, 0){4}}
\put(9, 1.7){\line(0, -1){1.7}}

\put(9, 2.5){\oval(0.8, 2.5)}
\put(6.9, 2){\oval(0.8, 4.5)}
\put(11.1, 2){\oval(0.8, 4.5)}

\qbezier(7, 0)(9, -1)(11, 0)
\end{picture}
\end{center}

\setlength{\unitlength}{0.59cm}
\begin{center}
\begin{picture}(5,7)

\put(0, 0){\circle*{0.15}}
\put(-0.6, 0){\footnotesize{$y_{4}$}}

\put(0, 2){\circle*{0.15}}
\put(-0.6, 2){\footnotesize{$y_{5}$}}

\put(0, 4){\circle*{0.15}}
\put(-0.6, 4){\footnotesize{$y_{6}$}}

\put(4, 0){\circle*{0.15}}
\put(4.2, 0){\footnotesize{$y_{3}$}}

\put(4, 2){\circle*{0.15}}
\put(4.2, 2){\footnotesize{$y_{2}$}}

\put(4, 4){\circle*{0.15}}
\put(4.2, 4){\footnotesize{$y_{1}$}}

\put(2, 1){\circle*{0.15}}
\put(1.9, 1.2){\footnotesize{$w'$}}

\put(2, 3){\circle*{0.15}}
\put(2, 2.5){\footnotesize{$w$}}

\put(0, 0){\line(0, 1){4}}
\put(0, 0){\line(1, 0){4}}
\put(4, 4){\line(-1, 0){4}}
\put(4, 4){\line(0, -1){4}}

\put(2, 2){\oval(1, 2.8)}

\put(0, 0){\line(2, 1){2}}
\put(4, 0){\line(-2, 1){2}}

\put(0, 4){\line(2, -1){2}}
\put(4, 4){\line(-2, -1){2}}

\put(-2, -1.2){\small{\textbf{Figure 4(c) :} The Class $\mathcal{F}_{3}$}}
\end{picture}
\end{center}
\vskip 30 pt

\noindent Let $P = u_{1}u_{2}u_{3}, P' = v_{1}v_{2}v_{3}$ and $P'' = w_{1}w_{2}w_{3}$ be three paths of length two. The graph $P_{3, 3, 3}$ is constructed from $P, P'$ and $P''$ by adding edges so that $\{u_{1}, v_{1}, w_{1}\}$ and $\{u_{3}, v_{3}, w_{3}\}$ form two complete graphs of order three. Brousek et al.~\cite{BRF} proved :

\begin{thm}\label{thm pumm}\cite{BRF}
Let $G$ be a $2$-connected $\{K_{1, 3}, N_{1, 2, 2}, N_{1, 1, 3}\}$-free graph. Then either $G$ is Hamiltonian, or $G$ is isomorphic to $P_{3, 3, 3}$ or $cl(G) \in \mathcal{F}_{1} \cup \mathcal{F}_{2} \cup \mathcal{F}_{3}$.
\end{thm}
\vskip 5 pt

\noindent Xiong et al.\cite{XLMWZ} generalized Theorem \ref{thm pumm} into $3$-connected graphs by establishing the following theorem.
\vskip 5 pt

\begin{thm}\cite{XLMWZ}\label{thm pum}
Let $G$ be a $3$-connected $\{K_{1, 3}, N_{s_{1}, s_{2}, s_{3}}\}$-free graph. If $s_{1} + s_{2} + s_{3} \leq 9$ and $s_{i} \geq 1$, then $G$ is Hamiltonian.
\end{thm}
\vskip 5 pt

\indent  In the context of $k$-$\gamma_{c}$-vertex critical graphs, Ananchuen et al. \cite{AAP1} established a basic property of these graphs.
\vskip 5 pt

\begin{lem}\label{lem 2}\cite{AAP1}
For $k \geq 2$, if $G$ is a $k$-$\gamma_{c}$-vertex critical graph and $v \in V(G)$, then $\gamma_{c}(G - v) = k - 1$.
\end{lem}
\vskip 5 pt

\noindent They, further, established a class of $3$-$\gamma_{c}$-vertex critical graphs. These graphs will be shown in this paper that they are $3$-$\gamma_{c}$-vertex critical non-Hamiltonian graphs containing a claw as an induced subgraph. For an integer $l \geq 6$, let $H_{l}$ be a copy of complete graph with the vertices $v_{1}, v_{2}, ..., v_{l}$ and remove a Hamiltonian cycle $v_{1}v_{2}...v_{l}$. Moreover, let $J$ be a graph of $\frac{l(l - 3)}{2} - l$ isolated vertices $v_{i, j}$ where $1 \leq i < j \leq l$ except the cases when $j = i + 1$ and $j = i + 2$ modulo $l$. The graph $J_{l}$ is constructed by $H_{l}$ and $J$ by adding edges from $H_{l}$ to $J$ except the edges $v_{i}v_{i, j}$ and $v_{j}v_{i, j}$.
\vskip 5 pt

\indent Recently, Kaemawichanurat and Caccetta~\cite{KC} established relationship between cardinalities of independent set and connected dominating set of claw free graphs.
\vskip 5 pt

\begin{lem}\label{lem P}\cite{KC}
Let $G$ be a claw-free graph, $X$ be a vertex subset and $I$ be an independent set of $G$ such that $X \succ_{c} I$. Then $|I| \leq |X| + 1$.
\end{lem}
\vskip 5 pt

\noindent In~\cite{KC1}, the author found a $5$-$\gamma_{c}$-vertex critical graph as illustrated by Figure 5. It is not difficult to see that this graph contains a claw as an induced subgraph and is non-Hamiltonian.

\setlength{\unitlength}{1cm}
\begin{center}
\begin{picture}(5, 5.3)
\put(0, 3){\circle*{0.15}}
\put(1.5, 1.5){\circle*{0.15}}
\put(1.5, 3){\circle*{0.15}}
\put(1.5, 4.5){\circle*{0.15}}
\put(3.5, 1.5){\circle*{0.15}}
\put(3.5, 2.5){\circle*{0.15}}
\put(3.5, 3.5){\circle*{0.15}}
\put(3.5, 4.5){\circle*{0.15}}
\put(4.5, 3){\circle*{0.15}}

\put(0, 3){\line(1, 0){1.5}}
\put(0, 3){\line(1, 1){1.5}}
\put(0, 3){\line(1, -1){1.5}}
\put(1.5, 1.5){\line(1, 0){2}}
\put(1.5, 4.5){\line(1, 0){2}}
\put(1.5, 3){\line(4, 1){2}}
\put(1.5, 3){\line(4, -1){2}}
\put(3.5, 1.5){\line(0, 1){3}}
\put(4.5, 3){\line(-2, 3){1}}
\put(4.5, 3){\line(-2, -3){1}}

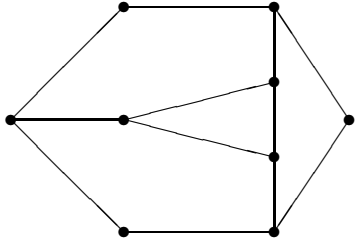
\put(-3, 0.3){\footnotesize{\textbf{Figure 5 :} A $5$-$\gamma_{c}$-vertex critical graph containing claw which is non-Hamiltonian.}}

\end{picture}
\end{center}
\vskip 15 pt

\indent Next, we will give a construction of some $k$-$\gamma_{t}$-vertex critical graphs from Wang and Wang~\cite{WW}. These graphs will be shown in this paper that they are $4$-$\gamma_{c}$-vertex critical non-Hamiltonian graphs containing a claw as an induced subgraph. Let
\vskip 5 pt

\indent $U_{1} = \{u_{1, 1}, u_{1, 2}, ..., u_{1, l}\}$, $U_{2} = \{u_{2, 1}, u_{2, 2}, ..., u_{2, l}\}$, $U_{3} = \{u_{3, 1}, u_{3, 2}, ..., u_{3, l}\}$, $U_{4} = \{y_{1}, y_{2}, y_{3}, y_{4}\}$.
\vskip 5 pt

\noindent A graph $G$ in the class $\mathcal{T}_{l}$ has the vertex set $U_{1} \cup U_{2} \cup U_{3} \cup U_{4} \cup \{u\}$ and the edge set
\vskip 5 pt

\indent $\{uu_{i, j} : i = 1, 2, 3, j = 1, 2, ..., l\} \cup \{y_{1}u_{1, i} : i = 1, 2, ..., l\}$
\vskip 5 pt

\indent $\cup \{y_{2}u_{2, i} : u_{2, i}, i = 1, 2, ..., l\} \cup \{y_{3}u_{i, j} : i = 2, 3, j = 1, 2, ..., k\}$
\vskip 5 pt

\indent $\cup \{y_{4}u_{3, i} : i = 1, 2, ..., l\} \cup \{y_{1}y_{2}, y_{2}y_{3}, y_{3}y_{4}, y_{4}y_{1}\}$
\vskip 5 pt

\indent $\{u_{1, i}u_{2, j} : i \neq j$ and $i, j = 1, 2, ..., l\} \cup \{u_{1, i}u_{3, j} : i \neq j$ and $i, j = 1, 2, ..., l\}$.
\vskip 5 pt

\indent We conclude this section with a result concerning the relationship between $4$-$\gamma_{c}$-vertex critical graphs and $4$-$\gamma_{t}$-vertex critical graphs. This result was proved by Kaemawichanurat et al.~\cite{KCA}.
\vskip 5 pt

\begin{thm}\label{thm mike}\cite{KCA}
Let $G$ be a $2$-connected graph. Then $G$ is $4$-$\gamma_{c}$-vertex critical if and only if $G$ is $4$-$\gamma_{t}$-vertex critical.
\end{thm}
\vskip 5 pt

\subsection{\bf Main results}
We are ready to prove the following theorems. For a graph $G$ containing a vertex $v$, we let $D^{c}_{v}$ be a smallest connected dominating set of $G - v$.
\vskip 5 pt

\begin{thm}\label{thm M}
Let $G$ be a $2$-connected $3$-$\gamma_{c}$-vertex critical claw-free graph. Then $G$ is Hamiltonian.
\end{thm}
\proof
We first show that $G$ is $\{N_{1, 2, 2}, N_{1, 1, 3}\}$-free graph. Suppose to the contrary that $G$ contains $N_{1, 2, 2}$ as an induced subgraph. Consider $G - w_{1}$. Clearly, the graph $G - w_{1}$ is claw free, moreover, $I_{1} = \{w_{2}, v_{1}, v_{3}, u_{1}\}$ is an independent set of $G - w_{1}$. Thus $|D^{c}_{w_{1}}| \succ_{c} I_{1}$. This implies by Lemma \ref{lem P} that $|D^{c}_{w_{1}}| + 1 \geq |I_{1}| \geq 4$. So $|D^{c}_{w_{1}}| \geq 3$ contradicting Lemma \ref{lem 2}. Therefore, $G$ is an $N_{1, 2, 2}$-free graph. Similarly, suppose to the contrary that $G$ contains $N_{1, 1, 3}$ as an induced subgraph. Consider $G - u_{1}$. Thus, the graph $G - u_{1}$ is claw free and $I_{2} = \{u_{2}, v_{1}, w_{1}, w_{3}\}$ is an independent set of $G - u_{1}$. So $|D^{c}_{u_{1}}| \succ_{c} I_{2}$. Lemma \ref{lem P} then gives that $|D^{c}_{u_{1}}| + 1 \geq |I_{2}| \geq 4$. Thus $|D^{c}_{w_{1}}| \geq 3$ contradicting Lemma \ref{lem 2}. Therefore, $G$ is an $N_{1, 1, 3}$-free graph.
\vskip 5 pt

\indent Clearly, $\gamma_{c}(P_{3, 3, 3}) = 5$. Therefore, $G$ is not $P_{3, 3, 3}$. To show that $G$ is Hamiltonian, by using Theorem \ref{thm pumm}, it remains to show that $cl(G) \notin \mathcal{F}_{1} \cup \mathcal{F}_{2} \cup \mathcal{F}_{3}$. Suppose to the contrary that $cl(G) \in \mathcal{F}_{1} \cup \mathcal{F}_{2} \cup \mathcal{F}_{3}$. We first consider the case when $cl(G) \in \mathcal{F}_{1}$. Consider $G - q_{1}$. Thus $|D^{c}_{q_{1}}| = 2$ by Lemma \ref{lem 2}. To dominate $V(Q_{2}) \cup V(Q_{3})$, we have $D^{c}_{q_{1}} \subseteq V(Q_{2}) \cup V(Q_{3})$. So $D^{c}_{q_{1}}$ does not dominate $V(Q_{1}) - \{q_{1}, z_{1}\}$, a contradiction. Therefore, $cl(G) \notin \mathcal{F}_{1}$. We now consider the case when $cl(G) \in \mathcal{F}_{2}$. Consider $G - c_{3}$. Similarly, $D^{c}_{c_{3}} \subseteq V(R_{1}) \cup V(R_{2})$. We see that $D^{c}_{c_{3}}$ does not dominate $r$. Thus $cl(G) \notin \mathcal{F}_{2}$. We finally consider the case when $cl(G) \in \mathcal{F}_{3}$. Consider $G - w$. Thus, $D^{c}_{w} \subseteq \{y_{1}, y_{2}, ..., y_{6}\}$ to dominate $\{y_{2}, y_{5}\}$. We see that $D^{c}_{w}$ does not dominate $K - \{w, w'\}$, a contradiction. Therefore, $cl(G) \notin \mathcal{F}_{3}$. In view of Theorem \ref{thm pumm}, $G$ is Hamiltonian.
\qed
\vskip 5 pt

\indent Recall the graph $J_{l}$ which is defined by Ananchuen et al.~\cite{AAP1} in Subsection \ref{SS2}. Observe that when $l \geq 8$, $|V(J)| > |V(H_{l})|$. Clearly, $V(H_{l})$ is a cut set of $J_{l}$ such that $J_{l} - V(H_{l})$ has the isolated vertices of $J$ as the $|J|$ components. Thus $\frac{|V(H_{l})|}{\omega(J_{l} - V(H_{l}))} < 1$. By Lemma \ref{lem Ch}, $J_{l}$ is non-Hamiltonian when $l \geq 8$. Therefore, the condition claw-free in Theorem \ref{thm M} is necessary.
\vskip 5 pt

\begin{thm}\label{thm W}
For $4 \leq k \leq 5$, let $G$ be a $3$-connected $k$-$\gamma_{c}$-vertex critical claw-free graph. Then $G$ is Hamiltonian.
\end{thm}
\proof
We will show that $G$ is an $N_{3, 3, 3}$-free graph. Suppose to the contrary that $G$ contains $N_{3, 3, 3}$ as an induced subgraph. Consider $G - w_{2}$. Clearly, the graph $G - w_{2}$ is claw-free and $I_{3} = \{w_{1}, w_{3}, u_{1}, u_{3}$ $, v_{1}, v_{3}\}$ is an independent set. Since $D^{c}_{w_{2}} \succ_{c} I_{3}$, by Lemma \ref{lem P}, $|D^{c}_{w_{2}}| \geq 5$. This implies by Lemma \ref{lem 2} that $k \geq 6 $, a contradiction. Hence, $G$ is an $N_{3, 3, 3}$-free graph. In view of Theorem \ref{thm pum}, $G$ is Hamiltonian.
\qed
\vskip 5 pt

\indent Recall the graph $\mathcal{T}_{l}$ which is defined by Wang and Wang~\cite{WW}. in Subsection \ref{SS2}. Theorem \ref{thm mike} implies that $\mathcal{T}_{l}$ is also a class of $4$-$\gamma_{c}$-vertex critical graphs. For a graph $G$ in this class, $G$ contains a claw as an induced subgraph and $U_{1} \cup U_{4} \cup \{u\}$ is a cut set of $l + 5$ vertices such that $G - (U_{1} \cup U_{4} \cup \{u\})$ has isolated vertices in $U_{2} \cup U_{3}$ as the $2l$ components. Thus, when $l \geq 6$, we have that $\frac{|U_{1} \cup U_{4} \cup \{u\}|}{\omega(G - (U_{1} \cup U_{4} \cup \{u\}))} < 1$. By Lemma \ref{lem Ch}, $G$ is non-Hamiltonian. Hence, the condition claw-free in Theorem \ref{thm W} is necessary.
\vskip 5 pt

\indent Finally, The graph in Figure 5 is a $5$-$\gamma_{c}$-vertex critical non-Hamiltonian graph containing a claw. Therefore, the condition claw-free and $3$-connected in Theorem \ref{thm W} are necessary.
\vskip 5 pt

\noindent \textbf{Acknowledgement:} The author thanks to assistant professor Thiradet Jiarasuksakun for his valuable advice. The author also express sincere thanks to the reviewer for greatly improving this paper.


\footnotesize
\begin{thebibliography}{99}\label{bib}
\bibitem{AP} N. Ananchuen and M. D. Plummer, Matching in $3$-vertex-critical: The odd case, Discrete Mathematics 307(2007) 1651 - 1658.
\bibitem{AAP1} W. Ananchuen, N. Ananchuen and M. D. Plummer, Vertex criticality for connected domination, Utilitas Mathematica 86(2011) 45 - 64.
\bibitem{BCD} R. C. Brigham, P. Z. Chinn and R. D. Dutton, Vertex domination-critical graphs, Networks 18(1988) 173 - 179.
\bibitem{B2} J. Brousek, Minimal $2$-connected non-hamiltonian claw-free graphs, Discrete Mathematics 191 (1998) 57 - 64.
\bibitem{BRF} J. Brousek, Z. Ryj\'{a}$\check{c}$ek and O. Favaron, Forbidden subgraphs, hamiltonicity and closure in claw-free graphs, Discrete Mathematics 196 (1999) 29 - 50.
\bibitem{C} V. Chv\'{a}tal, Tough graphs and hamiltonian circuits, Discrete Mathematics 306(2006) 910 - 917(reprinted from Discrete Mathematics 5(1973) 215 - 228).
\bibitem{FTZ} O. Favaron, F. Tian and Lei. Zhang, Independence and hamiltonicity in 3-dominaion-critical graphs, J Graph Theory 25(1997) 173 - 184.
\bibitem{FTWZ} E. Flandrin, F. Tian, B. Wei and L. Zhang, Some properties of 3-domination-critical graphs, Discrete Mathematics 205(1999) 65 - 76.
\bibitem{FHM} J. Fulman, D. Hanson and G. MacGillivray, Vertex domination-critical graphs, Networks 25(1995) 41 - 43.
\bibitem {HHS} T.W. Haynes, S.T. Hedetniemi, and P.J. Slater, Domination in graphs advanced topics, Marcel Dekker, Inc., New York, 1998.
\bibitem{R} Z. Ryj\'{a}$\check{c}$ek, On a closure concept in claw-free graphs, Journal of Combinatorial Theory(B) 70 (1997) 217 - 224.
\bibitem{KC1} P. Kaemawichanurat, Connected domination critical graphs. Ph.D. Thesis. PhD supervisor: Louis Caccetta. Curtin University (2016).
\bibitem{KC} P. Kaemawichanurat and L. Caccetta, Hamiltonicity of domination critical claw-free graphs, accepted by Journal of Combinatorial Mathematics and Combinatorial Computing.
\bibitem{KCA} P. Kaemawichanurat, L. Caccetta and N. Ananchuen, Critical graphs with respect to total domination and connected domniation, Australasian Journal of Combinatorics 65(1)(2016) 1 - 13.
\bibitem{KCA2} P. Kaemawichanurat, L. Caccetta and W. Ananchuen, Hamiltonicity of connected domination critical graphs. Ars Combinatoria. 136(2018) 127 - 151.

\bibitem{SB}D. P. Sumner and P. Blitch, Domination critical graphs, Journal of Combinatirial Theory(B) 34(1983) 65 - 76.
\bibitem{FWZ} F. Tian, B. Wei and L. Zhang, Hamiltonicity in 3-domination critical graphs with $\alpha = \delta + 2$, Discrete Applied Mathematics 92(1999) 57 - 70.
\bibitem{WW} Y. Wang and C. Wang, On $4$-$\gamma_{t}$-critical graphs of diameter two, Discrete Applied Mathematics 161(2013) 1660 - 1668.
\bibitem{W} E. Wojcicka, Hamiltonian properties of domination critical graphs, Journal of Graph Theory 14(1990) 205 - 215.
\bibitem{XLMWZ} W. Xiong, H. J. Lai, X. Ma, K. Wang and M. Zhang, Hamilton cycles in $3$-connected claw-free and net-free Graphs, Discrete Mathematics 313 (2013) 784 - 795.
\bibitem{YCXYX} Y. Yuansheng, Z. Chengye, L. Xiaohui, J. Yongsong and H. Xin, Some $3$-connected $4$-edge critical non-Hamiltonian graphs, Journal of Graph Theory 50: 316 - 320, 2005.
\end{thebibliography}
\end{document}